\documentclass[]{amsart}
\usepackage[]{amsmath, amsthm, amsfonts, amssymb }

\newcommand {\R} {{\mathbb R}}
\newcommand {\Z} {{\mathbb Z}}
\newcommand {\Q} {{\mathbb Q}}

\newtheorem{thm}{Theorem}
\newtheorem{cor}{Corollary}
\newtheorem{lemma}{Lemma}

\newtheorem{remark}{Remark}

\newtheorem{con}{Construction} 

\begin{document}

\title{ K\"ahler Solvmanifolds}

\author{
        Donu Arapura    
}

\address{Department of Mathematics\\
Purdue University\\
West Lafayette, IN 47907\\
U.S.A.}
\thanks{Partially supported by the NSF}
\email{arapura@math.purdue.edu}

\begin{abstract}
Compact K\"ahler solvmanifolds are classified up to biholomorphism.
A proof   of a  conjecture Benson and Gordon, that completely solvable
 compact K\"ahler solvmanifolds are tori is deduced from this.
\end{abstract}

\maketitle

In this note, we use the term
solvmanifold in a somewhat restricted sense.
A solvmanifold will be a compact quotient  $G/\Gamma$
of a connected simply connected solvable real Lie group
by a discrete subgroup.  Such a space is a $K(\Gamma,1)$.
Our goal is  to classify all K\"ahler solvmanifolds   up to biholomorphism.
To be precise, we consider a K\"ahler solvmanifold to be a
solvmanifold with   complex structure 
and  K\"ahler metric that are  not assumed invariant
under the group. (The invariance of the complex structure will follow
from our theorem.) Complex tori are clearly examples
of such spaces.  Hasegawa \cite{h2} has shown how to construct additional
examples by taking quotients of suitable tori, and we show that all
K\"ahler solvmanifolds  arise this way.
As a corollary, we prove a conjecture of Benson
and Gordon \cite{bg2} that tori are the only  K\"ahler
solvmanifolds corresponding to completely solvable Lie
groups\footnote{
The proof in \cite{tk} is incomplete; see \cite[p.257]{h2}.}.
The main new ingredient here is the restriction theorem for
polycyclic  K\"ahler groups proved by Nori and the author \cite{an}.

After the completion of the first version of this paper, I have
learned from Alex Brudnyi that he had obtained similar results in
his 1995 PhD thesis (Technion) and subsesquent papers \cite{br1, br2}.
However the methods are quite different and do not yield the 
analogous result for  solvmanifolds which are bimeromorphic to a compact K\"ahler
manifolds (theorem 2).

Let us start by defining a few basic terms and constructions.
A solvable group $G$ is called completely
solvable, or of real type,  if all the eigenvalues of images of the 
elements of its Lie algebra under its adjoint representation are real.
Given a solvmanifold $G/\Gamma$, the map 
$$G/\Gamma\to \frac{G/[G,G]}{im(\Gamma\to G/[G,G])}$$ 
defines a fiber bundle over a torus
called the canonical torus fibration in \cite{as}.

\begin{lemma}\label{lemma:torusfib}
  The fibers are compact nilmanifolds (i.e. quotients of connected nilpotent
Lie groups by lattices).
\end{lemma}

\begin{proof}
Since $[G,G]$ is nilpotent \cite[3.11]{r}, it lies in the 
maximum connected nilpotent subgroup $N$. $N\cap \Gamma$ is a
cocompact lattice in $N$ by a theorem of Mostow
\cite{m}, see also \cite[3.3]{r}. 
Hence  $[G,G]\cap \Gamma$ is also  a cocompact lattice in $[G,G]$
\cite[2.3]{r}.
\end{proof}

A complex torus is a quotient $V/L$ of a complex vector space by a lattice.
We can identify $V$ with the holomorphic tangent space at $0$, and the map
$V\to V/L$ with the exponential map. Holomorphic maps of complex tori
can always be lifted to affine maps of vector spaces \cite[p.326]{gh}; 
in particular, maps preserving $0$ are necessarily homomorphisms.
For any compact K\"ahler manifold $X$, the Albanese torus $Alb(X)$ is a complex
torus with real dimension equal to the first Betti 
number of $X$ \cite[p. 331]{gh}. After choosing a base point $x\in X$,
there is a holomorphic map, called
the Abel-Jacobi map,
$\alpha:X\to Alb(X)$ with $\alpha(x)=0$  which is topologically 
 just the classifying map $X\to K(H_1(X,\Z)/torsion,1)$.

\begin{lemma}\label{lemma:torus}
  If a compact K\"ahler manifold is diffeomorphic to a torus
then it is biholomorphic to a complex torus.
\end{lemma}

\begin{proof}
 Let $X$ be the manifold in question. Then it follows
 that the Albanese map is a holomorphic diffeomorphism
and hence biholomorphism. 
\end{proof}

When $X$ is a K\"ahler solvmanifold, the Abel-Jacobi map would
coincide with the canonical torus fibration. We will make use
of this fact later on.

We give  a  general method for producing K\"ahler
solvmanifolds along the lines of Hasegawa's paper \cite{h2}.
This method will produce all examples.
We will make use of the standard description of group extensions
using cohomology \cite{b}.

\begin{con}
Let
$$0\to \Z^{2m}\to \Gamma\to \Z^{2n}\to 0$$
be a group extension such that 
\begin{enumerate}
\item[(a)] The homomorphism $\Z^{2n}\to GL(2m,\Z)$, associated to the
extension, has finite image.
\item[(b)]  $\Z^{2n}\to GL(2m,\Z)$
extends to a homomorphism $\R^{2n}\to  GL(2m,\R)$ of
real Lie groups.
\item[(c)] There is a complex structure $J_0$ on $\R^{2m}$ invariant
 under the action of  $\Z^{2n}$.
\item[(d)] The class of the extension in $H^2( \Z^{2n},  \Z^{2m})$ has finite
  order.
\end{enumerate}
Let $V  =\R^{2m}\oplus \R^{2n}$, and choose complex structure
$J$ on it extending $J_0$.
Let $\Gamma\otimes \R$ denote the semidirect product  $\R^{2m}\rtimes \R^{2n}$
under the action (b). Note that as a sets,  $\Gamma\otimes \R$  and $V$ are
one and the same, nevertheless it is useful to keep a distinction.
Consider the diagram
$$
\left.
\begin{array}{ccccccccc}
 0 & \to & \Z^{2m} & \to & \Gamma & \to  & \Z^{2n} & \to & 0 \\
  &  & \downarrow &  & \downarrow &  & \downarrow &  &  \\
 0 & \to  & \R^{2m} & \to & \tilde\Gamma & \to & \Z^{2n} & \to & 0
\end{array}
\right.
$$
where the left hand square is the ``pushout''. Assumption
(d) shows that the bottom sequence splits. Therefore we get
embeddings $\Gamma\subset \tilde \Gamma\subset \Gamma\otimes \R$.
The quotient $X=(\Gamma\otimes \R)/\Gamma$ is a solvmanifold.
Under the identification $\Gamma\otimes \R= V$,
$\Gamma$ will act on $V$  by complex affine transformations.
Therefore $X$ inherits a complex structure. We can find,  thanks to (a),
a $\Gamma$-invariant Hermitean inner product on $V$, and this descends
to a flat K\"ahler metric on $X$. Note that different choices of $J$
can yield nonbiholomorphic, but diffeomorphic, examples.
\end{con}

Examples of the above type can be constructed quite
explicitly. For instance, they include all 
the  hyperelliptic (also called bielliptic) surfaces \cite{h2}.
In particular, they not all tori. However,  they are 
always quotients of tori by a finite group acting freely.
This is a consequence of the next lemma which implies that
$\Gamma$ contains a normal abelian subgroup $\Delta$ of finite
index. Then $X$ is a quotient of
the complex torus $V/\Delta$ by $\Gamma/\Delta$.

\begin{lemma}\label{lemma:ext}
Suppose that $\Gamma$ is a finitely generated group which 
fits into an extension, 
$$0\to \Z^{2m}\to \Gamma\stackrel{\pi}{\to} \Z^{2n}\to 0$$
Then conditions (a) and (d), above, hold if and only if $\Gamma$ contains a
normal abelian subgroup of finite index.
\end{lemma}

\begin{proof}
Suppose conditions (a) and (d) hold. Let $\Delta$ be the preimage in
$\Gamma$ of $\Delta'=ker[\Z^{2n}\to GL(2m,\Z)]$. This has finite
index. Since $H^2(\Delta',\Z^{2m})$ is torsion free \cite{b}, 
(d) shows that the sequence
$$0\to \Z^{2m}\to \Delta\to \Delta'\to 0$$
splits. Therefore $\Delta\subset\Gamma$ is a normal abelian subgroup of
finite index.

  Conversely, let $\Delta\subset \Gamma$ be a normal abelian subgroup of
  finite index. Since $\pi(\Delta) $ acts trivially on $\Delta\cap \Z^{2m}$
which has  finite index in $\Z^{2m}$, it follows that it
acts trivially on $\Z^{2m}$. Consequently  $\Z^{2n}\to GL(m,\Z)$
factors through $\Z^{2n}/\pi(\Delta)$. Let 
$\tilde\Delta=\pi^{-1}(\pi(\Delta))$ and consider the commutative
diagram
$$
\left.
\begin{array}{ccccccccc}
 0 & \to & \Z^{2m} & \to & \Gamma & \to & \Z^{2n} & \to & 0 \\
  &  & \| &  & \uparrow &  & \uparrow &  &  \\
 0 & \to & \Z^{2m} & \to & \tilde\Delta & \to & \pi(\Delta) & \to & 0 \\
  &  & \uparrow &  & \uparrow &  & \| &  &  \\
 0 & \to & \Delta\cap\Z^{2m} & \to & \Delta & \to & \pi(\Delta) & \to & 0
\end{array}
\right.
$$
Since $\Delta$ is abelian, the extension class of the
bottom row is $0$. The lower left hand square is a pushout.
Therefore the extension class  of the bottom row maps to the class
of the middle row, and consequently this is also zero.
The upper right hand square is a pullback. 
This implies that the extension class $e$ of the top row maps to the class
of the middle row under the map 
$H^2(\Z^{2n},\Z^{2m})\to H^2(\pi(\Delta),\Z^{2m}).$
This homomorphism is injective modulo torsion, since the transfer map \cite{b}
gives a splitting after tensoring with $\Q$. 
Therefore $e$ must be torsion.
\end{proof}

\begin{thm}\label{thm:1}
 Every  K\"ahler solvmanifold 
 is biholomorphic to one of the examples described in Construction 1.
\end{thm}    

\begin{remark} The theorem is a classification up  to biholomorphism
not isometry. These examples can have nonflat K\"ahler metrics,
and there does not appear to be any sensible way to classify these.
\end{remark}

\begin{proof}[Proof of theorem]
Let $\Gamma$ be a cocompact lattice in
a simply connected solvable Lie group $G$ such that $X=G/\Gamma$ has
a K\"ahler structure.
  By \cite[4.28]{r}, $\Gamma$ is polycyclic. Therefore \cite[4.10]{an}
implies that $\Gamma$ contains a nilpotent subgroup
$\Delta_1\subseteq\Gamma$ of finite index. We can  find a
subgroup $\Delta\subseteq \Delta_1$ of finite index which
is normal in $\Gamma$. $\Delta$ is again nilpotent, thus
$T=\Gamma/\Delta$ is an aspherical manifold with a nilpotent
fundamental group. Pulling back the K\"ahler structure from
$X$ also makes $T$ into a K\"ahler manifold such that the projection
$T\to X$ is holomorphic.
Therefore by \cite[theorem B]{bg1}, $T$ is homotopic to a torus.
Consequently, \cite[theorem A]{m} shows that $T$ is  diffeomorphic to 
a torus. Then $T$ is a complex torus by  lemma \ref{lemma:torus}.

 Consider the Abel-Jacobi  map $\alpha:X\to A = Alb(X)$.
This is surjective since it coincides with the canonical torus
 fibration. Let $f:T\to A$ denote the composition $T\to X\to A$. 
 We can choose a base point  of
 $X$ so that $0\in T$ maps to $0\in A$. Then $f:T\to A$ becomes a 
 homomorphism. Therefore $f$ is a bundle  with fiber
 $K = ker(f)$ which is a product of a torus with a finite abelian group.
 In particular, $T\to A$ is  a submersion,
 which implies that $X\to A$ is also a submersion. 
Therefore the fibers of $\alpha$ are connected K\"ahler
nilmanifolds (lemma~\ref{lemma:torusfib}), and hence tori by
\cite[theorem A]{bg1} or \cite{h1}.
It follows also that $[G,G]$ is abelian, and hence isomorphic
to a Euclidean space (with even dimension).

The exact sequence for homotopy groups of a fibration applied to 
$\alpha$ yields the extension
$$0\to \Z^{2m}\to \Gamma\to \Z^{2n}\to 0$$
We have to  check that the conditions of construction 1 hold.
Properties (a) and (d) follow from lemma \ref{lemma:ext}.
The action of of $G/[G,G]\cong \R^{2n}$ on $[G,G]\cong \R^{2m}$ 
gives the extension required for (b).
Let $V= \R^{2m}\oplus \R^{2n}$, $\Gamma$ acts on this by real affine
transformations. $V$ has a complex structure $J$ coming
from the identification of $V$ with the tangent space of  $T$ at $0$.
$\R^{2m}$ is invariant under $J$ since it can be identified with the 
tangent space to $K$ at $0$. The action of $\Gamma$ on $V$ is
holomorphic with respect to $J$, therefore $\Gamma$ acts by complex
affine transformations. In particular, the action of $\Z^{2n}$  on $\R^{2m}$
leaves $J|_{\R^{2m}}$ invariant, and this completes the proof.
\end{proof}

 As corollaries, we confirm  conjectures of Hasegawa and
 Benson and Gordon.

\begin{cor}
  A K\"ahler solvmanifold is both a quotient of a complex torus
and a bundle of complex tori over a complex torus.
\end{cor}

 \begin{proof}
  In the notation of the previous proof, $X$ is quotient of $T$ by
$\Gamma/\Delta$.  We have a submersion $\alpha:X\to A$.
  The map $T\to X$  induces a surjective homomorphism from
 $K$ to a fiber of $\alpha$ with a finite kernel $H$.  From this, we can deduce
 that $X\to A$ is a locally trivial bundle with fiber $K/H$.
 \end{proof}

\begin{cor}
  Let $X= G/\Gamma$ be a K\"ahler solvmanifold
with $G$ completely solvable.
Then $X$ is biholomorphic to a complex
torus. 
\end{cor}

\begin{proof}
The theorem shows that $G$ contains a lattice which is free abelian
of even rank. At  this point, we are in a position to argue as in
\cite{tk}.
 We can apply Saito's rigidity theorem \cite{s} to
conclude that $G$ is in fact abelian, and therefore that $X$
is a torus.
\end{proof}

We can prove a slightly stronger result by the same method.

\begin{thm}
  The examples given in Construction 1 are the only complex
  solvmanifolds which are bimeromorphic to a compact K\"ahler
manifolds.
\end{thm}

It is just a matter of going through each step and making
sure that everything still works. We indicate the main points:

\begin{itemize}
\item The key step that a nilmanifold   cannot be bimeromorphic to  
a compact K\"ahler manifolds unless it is a torus
is already supplied by Hasegawa \cite{h1}. 
 This is deduced by showing that 
a nilmanifold  cannot be formal unless it is a torus and then applying 
\cite{dgms}.
This argument yields a bit more, namely a generalization of  \cite[theorem
B]{bg1}, that tori are the only compact bimeromorphically K\"ahler
manifolds homotopic to a nilmanifold.

\item The above discussion about the Albanese carries over to the class
compact bimeromorphically K\"ahler manifolds without any modification.
See for example \cite[p. 337]{gpr}

\item Since fundamental groups of complex manifolds are invariant under
blow ups with smooth center, the classes of fundamental groups of
compact  bimeromorphically K\"ahler manifolds and compact  K\"ahler
manifolds coincide. Thus \cite[4.10]{an} also applies to the former.
\end{itemize}



\begin{thebibliography}{ABCD}
\bibitem[AN]{an} D. Arapura, M. Nori
{\em Solvable fundamental groups of algebraic varieties
and K\"ahler manifolds},
Compositio Math. (1999)

\bibitem[AS]{aa} L. Auslander, M. Auslander,
{\em Solvable Lie groups and   locally Euclidean Riemannian
spaces}, Proc. AMS (1958)  

\bibitem[AS]{as} L. Auslander, R. Szczarba
{\em Vector bundles over tori and noncompact solvmanifolds},
Amer. J. Math (1975)    
    
\bibitem[BG1]{bg1} C. Benson, C. Gordon,
{\em K\"ahler and symplectic nilmanifolds},
Topology (1988)

\bibitem[BG2]{bg2} C. Benson, C. Gordon,
{\em K\"ahler structures on compact solvmanifolds},
Proc. AMS (1990)

\bibitem[B]{b} K. Brown, {\em Cohomology of groups},
Springer-Verlag (1982)

\bibitem[Br1]{br1} A. Brudnyi, {\em Classification theorem for a class of
  flat connections and representations of Kähler groups},
  Michigan Math. J.  46  (1999)

\bibitem[Br2]{br2} A. Brudnyi, {\em Solvable representations of
K\"ahler groups}, preprint arXiv.math.CV/9911080

\bibitem[DGMS]{dgms} P. Deligne, P.  Griffiths, J. Morgan,
D. Sullivan, {\em Real homotopy theory of K\"ahler manifolds},
Invent. Math. 29 (1975) 

\bibitem[GH]{gh} P. Griffiths, J. Harris
{\em Principles of algebraic geometry},
Wiley (1978)

 \bibitem[GPR]{gpr} H. Grauert, T. Peternell, R. Remmert (eds.),
 {\em Several Complex Variables VI},
 Springer-Verlag (1991)

\bibitem[H1]{h1} K. Hasegawa
{\em Minimal models of nilmanifolds},
Proc. AMS (1989)

\bibitem[H2]{h2} K. Hasegawa
{\em A class of compact K\"ahlerian solvmanifolds and
a general conjecture},
Geom. Ded. (1999)


\bibitem[M]{m} G. Mostow
{\em Factor spaces of solvable Lie groups},
Ann. Math (1954)

\bibitem[R]{r} M. S. Raghunathan,
{\em Discrete subgroups of Lie groups},
Springer-Verlag (1972)

\bibitem[S]{s} M. Saito, 
{\em Sur certains groupes de Lie r\'esolubles. II},
Sci. Papers Coll. Gen. Ed. Univ. Tokyo 7 (1957) 

\bibitem[TK]{tk} A. Tralle, W. Kedra,
{\em Compact completely solvable K\"ahler solvmanifolds are tori},
Internat. Math. Res. Not. (1997)



\end{thebibliography}
\end{document}